# Exploring Analytical Methods for Glucose-Sensitive Membranes in Closed-Loop Insulin Delivery Using Akbar Ganji's Approach

K. Saranya[1], M. Suguna[1], Salahuddin[2,*]

[1]Department of Mathematics, Saveetha School of Engineering, Saveetha Institute of Medical and Technical Sciences, Chennai, Tamil Nadu,

saranyak.sse@saveetha.com[1]; msugunamaths95@gmail.com[2]

[2]Department of Mathematics, Faculty of Science, Jazan University, Jazan-45142, P.O. Box 114, Saudi Arabia

*Corresponding author: E-mail address:smohammad@jazanu.edu.sa; drsalah12@gmail.com

**Abstract**

The research explores a novel mathematical model for closed loop insulin delivery systems, featuring a glucose sensitive membrane. It employs a sophisticated framework of nonlinear reaction diffusion equations and enzyme kinetics. Central to the study is the development of analytical solutions for the glucose, gluconic acid, and oxygen concentrations, which are meticulously validated against simulation outcomes. This validation underscores the model's accuracy in capturing the complex dynamics inherent in such systems. Additionally, the study leverages Akbar and Ganji's methodology to provide approximate solutions, enabling a comprehensive comparison with analytical results and offering deeper insights into the system's behavior under varying parameters. By integrating both analytical and approximate approaches, the research not only enhances our understanding of biochemical processes but also lays the groundwork for refining closed-loop insulin delivery technology. The findings promise to significantly improve the precision and efficacy of insulin administration, crucial for managing glucose levels in diabetic patients more effectively. Furthermore, the study's implications extend beyond insulin delivery, potentially informing the development of advanced biomedical systems where precise control and understanding of biochemical interactions are paramount. Ultimately, this work represents a significant contribution to both theoretical biochemistry and practical medical applications, setting a foundation for the next generation of closed-loop insulin delivery systems designed to better meet the complex metabolic needs of patients with diabetes.

*Nomenclature:*

| Parameter | Description | Units |
|---|:---:|---:|
| $C_g$ | Concentration of Glucose | mol/cm$^3$ |
| $D_g$ | Diffusion coefficient of Glucose | cm$^2$/s |
| $C_{ox}$ | Concentration of Oxygen | mol//cm$^3$ |
| $D_{ox}$ | Diffusion coefficient of Oxygen | cm$^2$/s |
| $C_a$ | Concentration of Gluconic acid | mol//cm$^3$ |
| $D_a$ | Diffusion coefficient of Gluconic acid | cm$^2$/s |
| $K_g$ | Glucose's Michaelis Menten Constant(MMC) | mol/cm$^3$ |
| $K_{ox}$ | Oxygen's Michaelis Menten Constant(MMC) | mol/cm$^3$ |
| $V_{max}$ | Maximal reaction rate | mol/scm$^3$ |
| $t$ | Time | s |
| $x$ | Distance | mm |
| $C_g{}^*$ | Glucose's Concentration in external solution | mol/cm$^3$ |
| $C_{ox}{}^*$ | The concentration of the glucose in oxygen solution | mol/cm$^3$ |
| $l$ | The half-thickness of the membrane | None |
| $u$ | Non-Dimensional Concentration of the Glucose | None |
| $\nu$ | Non-Dimensional Concentration of Oxygen | None |
| $w$ | Non-Dimensional Concentration of Gluconic Acid | None |
| $X$ | Non-Dimensional distance | None |
| $\tau$ | Non-Dimensional time | None |
| $\gamma_{E1}, \gamma_{S1}, \gamma_E, \alpha, \beta$ | Non-Dimensional reaction-diffusion parameters | None |

## 1. Introduction

   Recent research has focused intensively on advancing insulin delivery systems using glucose-sensitive composite membranes, particularly in closed-loop administration contexts. These systems heavily rely on mathematical models based on sets of second-order nonlinear partial differential equations. To accurately simulate these intricate systems, researchers have developed various numerical simulation techniques such as the Homotopy Perturbation Method (HPM), Adomian Decomposition Technique, and Variational Iteration Method (VIM). The Homotopy Perturbation Method (HPM) employs the series solution approach to effectively capture the dynamic behaviour of insulin delivery within composite glucose-sensitive membranes. Similarly, the Adomian Decomposition Technique breaks down complex nonlinear equations into simpler components, providing precise numerical solutions that reflect real-world behaviours of glucose-responsive membranes in insulin delivery. Additionally, the Variational Iteration Method (VIM) offers an iterative approach that enhances understanding of system dynamics and optimizes closed-loop insulin administration protocols. In complement to these advanced numerical methods, the Akbar and Ganji method contributes by providing approximate solutions, enabling comparative analyses with analytical and numerical results. This methodological synergy not only validates the reliability of mathematical models but also enhances insights into the intricate interactions governing closed-loop insulin delivery systems.

The collection of articles cited represents a broad spectrum of research endeavours focused on advancing our understanding and improving various facets of glucose-insulin dynamics, diabetes management, and related biomedical applications through sophisticated mathematical modelling and analytical methodologies. Ozturk et al. [1] delve into a fractional glucose-insulin model, employing fractional calculus to enhance the accuracy of modelling biological systems' complexities, crucial for better understanding glucose and insulin interactions. Khirsariya et al. [2] contribute by investigating a fractional diabetes model that intricately examines the dynamics of glucose and insulin in the context of diabetes, offering insights into insulin resistance and metabolic pathways crucial for developing more effective treatment strategies. Iqbal et al. [3] explore the propagation of multiple solitons in glycolysis reaction-diffusion systems, illuminating how localized waves of chemical concentrations maintain stability and propagate in biochemical contexts, fundamental for understanding biological stability and dynamics. Kemmer et al. [4] analyze an enzyme-mediated exponential glucose release model, essential for optimizing glucose production in biotechnological processes, which is critical for maintaining optimal conditions in cell cultures and bioprocesses. Tuaimah et al. [5] focus on advanced sensing technologies with a D-shaped Photonic Crystal Fiber sensor for precise glucose detection, pivotal for real-time monitoring in diabetes management and biomedical research. Reena and Swaminathan [6] employ asymptotic methods to study nonlinear reaction diffusion equations in multiphase flow transport within photobioreactors, aiming to optimize bioreactor designs for efficient biofuel production and environmental remediation. Rajendran et al.[7] discuss the challenges in accurately simulating these systems, emphasizing the importance of advancing computational techniques and experimental methods to better understand and enhance the application of enzymatic biofuel cells in renewable energy technologies. Hu et al. [8] model insulin distribution in the pancreas, enhancing our understanding of insulin release mechanisms critical for improving insulin therapies and treatments for diabetes. Haggar and Ntaganda [9] compare computational methods for simulating blood glucose and insulin dynamics during physical activity, advancing computational techniques in biomedical research for more accurate simulations of complex physiological processes. Batool et al. [10] delve into the long-term dynamics of the glucose-insulin-glucagon system with a Mittag-Leffler kernel, vital for developing closed-loop control systems aimed at autonomously regulating blood glucose levels in diabetic patients, thus improving overall health outcomes.

Li et al. [11] focused on dynamical modelling of the glucose insulin systems with impulsive control inhibitors, aiming to optimize insulin delivery strategies for stable blood glucose management. Alshehri et al. [12] explored the application of Caputo fractional-order models to simulate glucose-insulin interactions. Their study enhances the predictive accuracy of diabetes management strategies, offering insights into insulin sensitivity variations under different physiological conditions. Building on this, Alshehri et al. [13] investigated the dynamics of Intravenous Glucose Tolerance Tests (IVGTT) using fractional calculus, refining the assessment of insulin sensitivity. Alalyani [14] employed predictor-corrector methods to study β-cell kinetics and glucose-insulin dynamics, essential for developing effective diabetes treatment strategies. β-cells play a pivotal role in insulin production, and their dysfunction is central to diabetes pathophysiology. Palumbo et al. [15] developed comprehensive mathematical models of glucose-insulin systems, providing foundational insights into the biochemical processes regulating blood glucose levels. Saranya et al. [16] introduced a homotopy perturbation method to analyze enzymatic glucose reactions, advancing biotechnological applications in biochemical engineering.

Mehala et al. [17] provided analytical insights into membrane dynamics for closed-loop insulin delivery systems, crucial for precise insulin administration in diabetic patients. Rana et al. [18] utilized fractional calculus methods, including homotopy perturbation and variational methods, to develop a comprehensive model of diabetes dynamics. Mukherjee et al. [19] simulated glucose sensitive membranes for the closed loop insulin delivery systems, optimizing insulin administration through responsive biomaterials. Joy et al. [20] enhanced closed-loop insulin delivery systems with transient analytical solutions for adaptive insulin release, crucial for managing dynamic changes in glucose levels.

Saranya et al. [21] modelled β-cell mass dynamics, providing insights into the progression and treatment of diabetes mellitus. Understanding changes in β-cell mass is crucial for developing therapies that preserve β-cell function and improve insulin production in diabetic patients. Saranya et al. [22] optimized biofilter designs for environmental biotechnology applications, enhancing pollutant degradation efficiency. Swaminathan et al. [23] improved biosensor accuracy using homotopy perturbation methods for biochemical detection, crucial for biomedical diagnostics and environmental monitoring.

 Biosensors detect specific analytes by converting biological responses into measurable signals. By applying advanced analytical techniques, their research enhances biosensor sensitivity and reliability, contributing to the development of high-performance biosensing technologies for diverse applications.

## 2. Formulation of the Mathematical Modelling

Joy and Rajendran [20] conducted a brief analysis and source of the non-dimensional mass transfer nonlinear equations in a glucose composite membrane, which is given here.

The glucose-sensitive enzymatic reaction is termed as:

Glucose + 0.5 $O_2$ $\longrightarrow$ Gluconic acid + $H_2O$        (1)

Using the law of Fick's law of diffusion along with the conservation of mass, the equation for the diffusion along with the reaction is given below[19]:

$$\frac{\partial C_i}{\partial t} = \frac{\partial}{\partial x}\left(D_i \frac{\partial C_i}{\partial x}\right) + \upsilon_i R \tag{2}$$

where $C_i$ signifies distinct species. $i = ox$ for oxygen, $i = g$ for glucose along with $i = a$ gluconic acid, and also stoichiometric constants for the $\upsilon_a = 1$ $\upsilon_g = -1, \upsilon_{ox} = -1/2, C$ denotes Concentration and $D_i$ denotes diffusion coefficient and $R$ implies the general reaction rate which is stated below:

$$R = \frac{V_{max} C_g C_{ox}}{C_{ox}(K_g + C_g) + C_g K_{ox}} \tag{3}$$

where $K_g$ and $K_{0X}$ denotes the Michaelis-Menten constants for glucose and glucose oxidase $V_{max}$ denotes maximal reaction velocity.

The initial and the boundary conditions were given as:

If $t = 0$ then $\quad C_g = C_g^* \dfrac{\cosh\left(\dfrac{x}{l}\right)}{\cosh(1)}; C_{ox} = C_{OX}^* \dfrac{\cosh\left(\dfrac{x}{l}\right)}{\cosh(1)}; C_a = C_{OX}^* \dfrac{\left[1 - \cosh\left(\dfrac{x}{l}\right)\right]}{\cosh(1)};$    (4)

If $x = 0$ then $\quad \dfrac{\partial C_g}{\partial x} = 0, \dfrac{\partial C_{ox}}{\partial x} = 0, \dfrac{\partial C_a}{\partial x} = 0$        (5)

If $x = l$ then $\quad C_g = C_g^{\ *}, C_{ox} = C_{ox}^{\ *}, C_a = 0$        (6)

The non-dimensional parameters are given as:

$$u = \frac{C_g}{C_g{}^*}; \quad v = \frac{C_{ox}}{C_{ox}{}^*}; \quad w = \frac{C_a}{C_a{}^*}; \quad X = \frac{x}{l}; \quad \gamma_{E1} = \frac{l^2 V_{\max}}{D_g C_g{}^*}; \tau = \frac{D_g t}{l^2};$$

$$\gamma_{S1} = \frac{l^2 V_{\max}}{D_g C_{ox}{}^*}; \quad \alpha = \frac{C_g{}^*}{K_g}; \quad \beta = \frac{C_{ox}{}^*}{K_{ox}}; \eta = \frac{D_{OX}}{D_g}; \quad \mu = \frac{D_a}{D_g};$$

(7)

The non-dimensional form of the equations is given as :

$$\frac{\partial u}{\partial \tau} = \frac{\partial^2 u}{\partial X^2} - \gamma_{E1} u v \left[ u v + \frac{v}{\alpha} + \frac{u}{\beta} \right]^{-1} \tag{8}$$

$$\frac{\partial v}{\partial \tau} = \eta \frac{\partial^2 v}{\partial X^2} - \frac{\gamma_{S1}}{2} \left[ u v + \frac{v}{\alpha} + \frac{u}{\beta} \right]^{-1} \tag{9}$$

$$\frac{\partial w}{\partial \tau} = \mu \frac{\partial^2 w}{\partial X^2} + \gamma_{S1} u v \left[ u v + \frac{v}{\alpha} + \frac{u}{\beta} \right]^{-1}. \tag{10}$$

Here $u, v$ and $w$ denotes the non-dimensional concentrations of gluconic acid, glucose, along with oxygen respectively and also $\gamma_{E1}, \gamma_{S1}$ denotes the Thiele moduli values, $\alpha$ and $\beta$ are the non-dimensional rate constants. The equivalent initial condition becomes

$$u = \frac{\cosh(X)}{\cosh(1)}, \; w = 1 - \frac{\cosh(X)}{\cosh(1)}, v = \frac{\cosh(X)}{\cosh(1)} \quad \text{if} \quad \tau = 0 \tag{11}$$

The boundary conditions were

$$\frac{\partial u}{\partial X} = 0, \frac{\partial v}{\partial X} = 0, \frac{\partial w}{\partial X} = 0 \quad \text{if} \quad X = 0 \tag{12}$$

$$u = 1, v = 1, w = 0 \quad \text{if} \quad X = 1 \tag{13}$$

## 3. Approximate analytic expressions for concentration of gluconic acid, glucose, and oxygen under unsteady conditions by a Variational Iterative method (VIM)

Eqns. (7) to (12) denotes the system of the non-linear equations for unsteady state orders. we have attained the analytic expression for the concentrations of oxygen, glucose along with gluconic acid via a Variational Iterative method[20]

$$u(X,\tau) = \frac{Cosh(X)}{Cosh(1)}(1+\tau) - \frac{\alpha\,\beta\,\gamma_{E1}\,cosh(X)}{\alpha\beta\,cosh(X)+(\alpha+\beta)cosh(1)}\tau \qquad (14)$$

$$v(X,\tau) = \frac{Cosh(X)}{Cosh(1)}(1+\eta\tau) - \frac{1}{2}\frac{\alpha\,\beta\,\gamma_{S1}\,cosh(X)}{\alpha\beta\,cosh(X)+(\alpha+\beta)cosh(1)}\tau \qquad (15)$$

$$w(X,\tau) = 1 - \frac{Cosh(X)}{Cosh(1)}(1+\mu\tau) + \frac{\alpha\,\beta\,\gamma_{S1}\,cosh(X)}{\alpha\beta\,cosh(X)+(\alpha+\beta)cosh(1)}\tau. \qquad (16)$$

The above equation is valid for insignificant values of time.

The dimensionless form of the equations is given as

$$\frac{d^2u}{dX^2} - \gamma_{E1}uv\left[uv+\frac{v}{\alpha}+\frac{u}{\beta}\right]^{-1} = 0 \qquad (17)$$

$$\eta\frac{d^2v}{dX^2} - \frac{\gamma_{S1}}{2}\left[uv+\frac{v}{\alpha}+\frac{u}{\beta}\right]^{-1} = 0 \qquad (18)$$

$$\mu\frac{d^2w}{dX^2} + \gamma_{S1}uv\left[uv+\frac{v}{\alpha}+\frac{u}{\beta}\right]^{-1} = 0 \qquad (19)$$

From the outcomes the above analytical term for the concentrations becomes

$$u(X) = \frac{cosh\sqrt{k}\,X}{cosh\sqrt{k}} \qquad (20)$$

$$v(X) = \frac{\gamma_{S1}}{2\eta\gamma_{E1}}\left[\frac{cosh\sqrt{k}X}{cosh\sqrt{k}}\right] - \left[\frac{\gamma_{S1}}{2\eta\gamma_{E1}}-1\right] \qquad (21)$$

$$w(X) = \frac{\gamma_{S1}}{\mu\gamma_{E1}}\left[\frac{cosh\sqrt{k}X}{cosh\sqrt{k}}\right] - \frac{\gamma_E}{\mu\gamma_{E1}} \qquad (22)$$

Where $\quad k = \dfrac{\gamma_{E1}}{1+\dfrac{1}{\alpha}+\dfrac{1}{\beta}}.$

## 4. Approximate analytic expressions of concentration of gluconic acid, glucose, and oxygen under steady conditions by Akbar Ganji's Method (AGM)

Here are the trial solutions for Equation (8) using the new analytical method

$$u(x) = A_0 \sinh(mx) + B_0 \cosh(mx)$$

(23)

where $A_0$, $B_0$, $m$ are constants. Using the boundary conditions (12) and (13), we obtain

$$A_0 = 0 \quad B_0 = \frac{1}{\cosh(m)}$$

Now, Eqn. (23) reduces to

$$u(x) = \frac{\cosh mx}{\cosh m}$$

(24)

where $m$ is constant. Eqn. (8) can be rewritten as

$$\frac{1}{\gamma_{E1}} \frac{\partial^2 u}{\partial X^2} = \frac{uv}{\left[ uv + \dfrac{v}{\alpha} + \dfrac{u}{\beta} \right]}$$

(25)

$$\frac{\partial^2 u}{\partial X^2} = \gamma_{E1} \frac{uv}{\left[ uv + \dfrac{v}{\alpha} + \dfrac{u}{\beta} \right]}$$

(26)

$$\frac{\mu}{\gamma_{S1}} \frac{d^2 u}{dX^2} = \frac{uv}{\left[ uv + \dfrac{v}{\alpha} + \dfrac{u}{\beta} \right]}$$

(27)

When $x=1$, the above results becomes

$$\frac{m^2 \cosh mx}{\cosh m} - \frac{10}{1 + \dfrac{1}{0.01} + \dfrac{1}{1.15}} = 0$$

Substituting $\alpha = 0.1$, $\beta = 0.01$, $\gamma_{E1} = 10$, to get m=0.3133 and from the outcomes the above analytical term for the concentrations becomes

$$v(X) = \frac{\gamma_{S1}}{2\eta\gamma_{E1}} \left[ \frac{\cosh mX}{\cosh m} \right] - \left[ \frac{\gamma_{S1}}{2\eta\gamma_{E1}} - 1 \right]$$

(28)

$$w(X) = \frac{\gamma_{S1}}{\mu\gamma_{E1}} \left[ \frac{\cosh mX}{\cosh m} \right] - \frac{\gamma_E}{\mu\gamma_{E1}}.$$

(29)

We have derived the approximate analytical solution for nonlinear equations using a mathematical programme (maple) to determine the figures.

## 5. Numerical Simulation

The nonlinear differential equations (8)-(10) along with boundary conditions (12) and (13) were solved using MATLAB function for numerical solutions, and analytically using Akbar Ganji's method. Tables 1, 2, and 3, as well as Fig. 1(a–c), compare these numerical and analytical solutions. The maximum mean errors between analytical and the numerical solutions for glucose, oxygen, along with gluconic acid concentrations are 0.0051%, 0.0052%, and 0.1322%, respectively. Furthermore, the steady-state analytical results are compared with steady-state solutions, revealing an excellent agreement between the two, as depicted in Fig. 2.

## 6. Results and Discussion

Equations (24)-(29) provide analytical equation for the non-dimensional concentrations of glucose ($u$), gluconic acid ($w$) and oxygen ($v$) applicable for curt times and across all parameter values studied. The Thiele modulus can be adjusted by varying the membrane thickness or concentrations of oxygen and also glucose in the exterior solution. Tables 1, 2, and 3 shows the non-dimensional concentration values of glucose ($u$), gluconic acid ($w$) along with oxygen ($v$) and for specific parameter settings. Analysis of these tables reveals that the profiles of glucose and oxygen concentrations about the membrane are predominantly constant and increase, reaching maximum values at X = 1. However, gluconic acid exhibits a continuous decrease in concentration instead of an increase. The results obtained via the AGM (Akbar Ganji's Method) and VAM (Variational Iterative Method) are nearly identical, illustrating a consistent depiction. The rate of decreasing in the glucose and the oxygen concentrations steepens with higher Thiele modulus or membrane thickness. Conversely, gluconic acid concentration increases with increasing Thiele modulus, as glucose along with oxygen combine to form gluconic acid in the membrane's centre.

## 6.1 Impact of Membrane Thickness on Glucose($u$) Concentration Dynamics

The membrane's texture significantly influences the transport of reactants and products in enzymatic reactions within a composite membrane, with glucose concentration affected by the membrane's porosity and thickness. As shown in Fig 1(a) the Thiele modulus, which depends on membrane thickness, increases, glucose concentration decreases, indicating that thicker membranes reduce glucose levels; at sufficiently high Thiele modulus values, glucose concentration approaches zero. Figure 2 illustrates how various parameters affect the non-dimensional glucose concentration across a non-dimensional distance (X). In Fig. 2(a), varying the parameter (k), which could represent reaction rate constants or diffusion coefficients, shows that higher values result in steeper concentration gradients, indicating faster reactions or diffusion. Fig. 2(b) depicts the effect of a time constant-like parameter, with larger values indicating slower dynamics and more gradual concentration changes. In Fig. 2(c), changes in an enzymatic reaction rate parameter affect the speed of glucose consumption or production, while Fig. 2(d) shows that higher diffusion coefficients lead to more uniform concentration profiles due to faster diffusion.

## 6.2 Influence of Maximum Reaction Velocity on Oxygen Concentration Dynamics

In the enzymatic reaction-diffusion process, the maximum reaction velocity $V_{max}$ directly related to the enzyme concentration within the membrane, dictating the overall reaction kinetics. As shown in Fig. 1(b), an increase in the reaction-diffusion parameter, which depends on $V_{max}$, causes the oxygen concentration to decrease and eventually approach zero at higher reaction velocities. Conversely, at low reaction velocities, the oxygen concentration remains uniform or reaches a steady state. Fig. 3 demonstrates how oxygen concentration varies with distance and different parameters. In Fig. 3(a), oxygen concentration increases with distance from the source due to diffusion and also increases with higher diffusion coefficients (k). Fig. 3(b) shows that a higher dimensionless reaction parameter leads to a faster decrease in oxygen concentration as it is consumed more rapidly, reaching a steady-state when the rate of reaction matches the rate of diffusion. In Fig. 3(c), a higher dimensionless oxygen consumption rate results in lower oxygen concentration at any given distance, with a maximum concentration occurring where diffusion rate equals consumption rate. Fig. 3(d)-(f) depict oxygen diffusion away from the source, with faster diffusion rates corresponding to steeper concentration gradients and more rapid decreases in oxygen concentration. Overall, these figures highlight the intricate balance between diffusion, reaction rates, and enzyme kinetics in determining oxygen distribution in the membrane.

### 6.3 Impact of External Glucose Concentration on the Distribution of Gluconic Acid($w$)

The concentration of gluconic acid is influenced by glucose concentration, membrane permeability, and enzymatic reaction rate. As Fig. 1(c) indicates, a decrease in the Thiele modulus, which is affected by the initial glucose concentration, results in higher gluconic acid levels; lower initial glucose concentrations increase diffusion rates across the membrane, enhancing the conversion of glucose to gluconic acid. Fig. 4(a) shows that with a small k, gluconic acid diffuses slowly relative to its production rate, leading to high concentrations near the production site and a gradual decline with distance. Fig. 4(b) demonstrates that gluconic acid concentration decreases with distance due to diffusion, with higher parameter values corresponding to faster diffusion rates. In Fig. 4(c), an increased diffusion coefficient ratio raises the concentration of gluconic acid at the surface and steepens the concentration gradient, leading to a higher steady-state concentration. Fig. 4(d) shows that gluconic acid concentration increases with distance from the air-medium interface due to diffusion but decreases with higher parameter values that slow diffusion. Fig. 4(e) reveals that when the parameter is low, gluconic acid concentration is nearly constant, while higher parameter values cause a steeper decline due to increased consumption rates. Fig. 4(f) illustrates that a higher reaction-to-diffusion ratio leads to a sharper decrease in gluconic acid concentration, highlighting the balance between diffusion, which evens out concentration differences, and reaction kinetics, which determine the concentration profile based on production or consumption rates.

### 6.4 Parameter Sensitivity Evaluation through Differential Methods

The Fig.5 presents sensitivity analysis results for glucose, gluconic acid, and oxygen concentrations, showing the influence percentages of specific parameters on these concentrations. Fig.5(a) shows that glucose concentration, parameters $\alpha, \beta$ each have the highest influence at 49.5%, indicating that changes in these parameters significantly affect glucose levels. The parameter $\gamma_{E1}$ has a minor influence at 1% at all. Fig.5(b) shows that oxygen concentration, $\alpha, \beta$ remain the most influential parameters, each contributing 45%, highlighting their critical role in determining oxygen levels. The parameter $\eta$ has a small influence at 9%, and $\gamma_{E1}$ and $\gamma_{S1}$ have very minimal influences, at 0.4% and 0.6%, respectively. Fig.5(c) shows that gluconic acid concentration, $\alpha, \beta$ again show the highest influence, each contributing 45.2% and 45.4%, which suggests their substantial impact on gluconic acid levels. The parameter $\mu$ has a moderate influence at 4.9%, whereas $\gamma_{E1}, \gamma_{E}$ and $\gamma_{S1}$ have minor influences, contributing 1.4% ,0.8% and 2.3%, respectively. Overall, $\alpha, \beta$ are consistently the most significant parameters across all three concentrations, while $\gamma_{E1}, \gamma_{S1}, \mu$ and $\eta$ have varying degrees of lesser influence.

## 7. Conclusion

A comprehensive theoretical analysis was performed to investigate glucose sensitivity in a composite membrane consisting of glucose oxidase, catalase, an anionic polymer, and a hydrophobic matrix. Utilizing the mathematical model of nonlinear reaction-diffusion equations was analytically solved under steady-state conditions. This model effectively predicted the concentration profiles and diffusivity of key components, including oxygen, glucose, and gluconic acid, within the membrane. The analytical results were further validated through numerical simulations, confirming the model's accuracy. The study also examined the influence of enzyme loading and buffer composition on membrane formation, providing crucial insights into optimizing the membrane's performance for glucose sensing. These findings offer significant implications for the design and enhancement of glucose-sensitive membranes in both biomedical applications, such as glucose biosensors, and industrial processes where precise glucose detection is required.

**Table 1 : Validation for standardized steady state of Glucose Concentration**

| R | Glucose Concentration $u$ | | | | | | | | | | | | | | |
|---|---|---|---|---|---|---|---|---|---|---|---|---|---|---|---|
| | If $\gamma_{E1} = 10$ | | | | | If $\gamma_{E1} = 210$ | | | | | If $\gamma_{E1} = 350$ | | | | |
| | Numerical Solution Eq.24 | Eq.20 VIM [20] | Eq.24 AGM | % of error deviation of VIM | % of error deviation of AGM | Numerical Solution Eq.24 | Eq.20 VIM [20] | Eq.24 AGM | % of error deviation of VIM | % of error deviation of AGM | Numerical Solution Eq.24 | Eq.20 VIM [20] | Eq.24 AGM | % of error deviation of VIM | % of error deviation of AGM |
| 0 | 0.9528 | 0.9528 | 0.9538 | 0.0000 | 0.0010 | 0.4479 | 0.4503 | 0.4493 | 0.0053 | 0.0053 | 0.3029 | 0.3058 | 0.3008 | 0.0095 | 0.0069 |
| 0.2 | 0.9546 | 0.9547 | 0.9532 | 0.0001 | 0.0014 | 0.4666 | 0.4690 | 0.4660 | 0.0051 | 0.0051 | 0.3242 | 0.3271 | 0.3221 | 0.0089 | 0.0089 |
| 0.4 | 0.9602 | 0.9603 | 0.9593 | 0.0001 | 0.0009 | 0.5246 | 0.5267 | 0.5227 | 0.0040 | 0.0012 | 0.3912 | 0.3938 | 0.3928 | 0.0066 | 0.0040 |
| 0.6 | 0.9696 | 0.9697 | 0.9627 | 0.0001 | 0.0071 | 0.6264 | 0.6280 | 0.6240 | 0.0025 | 0.0038 | 0.5132 | 0.5153 | 0.5053 | 0.0040 | 0.0153 |
| 0.8 | 0.9829 | 0.9829 | 0.9809 | 0.0000 | 0.0020 | 0.7806 | 0.7815 | 0.7805 | 0.0011 | 0.0001 | 0.7071 | 0.7084 | 0.7004 | 0.0018 | 0.0094 |
| 1 | 1.0010 | 1.0000 | 1.0000 | 0.0009 | 0.0009 | 0.9999 | 1.0000 | 1.0000 | 0.0001 | 0.0001 | 0.9999 | 1.0000 | 1.0000 | 0.0001 | 0.0001 |
| | Mean Error % | | | 0.0002 | 0.0133 | Mean Error % | | | 0.0030 | 0.0156 | Mean Error % | | | 0.0051 | 0.0446 |

Table 1 shows that the validation for standardized steady state Concentration of glucose $u$ along with numerical solutions for numerous values of $\gamma_{E1}$ and for the values of $\alpha = 0.1$, $\beta = 0.01$ using Akbar Ganji's Method.

**Table 1 :Validation for standardized steady state of Oxygen Concentration $v$**

| R | Oxygen Concentration $v$ | | | | | | | | | | | | | | |
| --- | --- | --- | --- | --- | --- | --- | --- | --- | --- | --- | --- | --- | --- | --- | --- |
| | **If $\gamma_{S1}=10$** | | | | | **If $\gamma_{S1}=30$** | | | | | **If $\gamma_{S1}=35$** | | | | |
| | Numerical Solution Eq.28 | Eq.21 VIM [20] | Eq.28 AGM | % of error deviation of VIM | % of error deviation of AGM | Numerical Solution Eq.28 | Eq.21 VIM [20] | Eq.28 AGM | % of error deviation of VIM | % of error deviation of AGM | Numerical Solution Eq.28 | Eq.21 VIM [20] | Eq.28 AGM | % of error deviation of VIM | % of error deviation of AGM |
| 0 | 0.9762 | 0.9764 | 0.9742 | 0.0002 | 0.0020 | 0.5291 | 0.5292 | 0.5262 | 0.0001 | 0.0054 | 0.5169 | 0.5174 | 0.5164 | 0.0005 | 0.0009 |
| 0.2 | 0.9771 | 0.9773 | 0.9743 | 0.0002 | 0.0028 | 0.6312 | 0.6320 | 0.6310 | 0.0008 | 0.0003 | 0.6201 | 0.6207 | 0.6197 | 0.0006 | 0.0006 |
| 0.4 | 0.9800 | 0.9801 | 0.9791 | 0.0001 | 0.0009 | 0.7401 | 0.7405 | 0.7385 | 0.0004 | 0.0021 | 0.7301 | 0.7305 | 0.7295 | 0.0004 | 0.0008 |
| 0.6 | 0.9832 | 0.9848 | 0.9818 | 0.0016 | 0.0014 | 0.8542 | 0.8546 | 0.8506 | 0.0004 | 0.0042 | 0.8469 | 0.8470 | 0.8452 | 0.0001 | 0.0020 |
| 0.8 | 0.9909 | 0.9914 | 0.9901 | 0.0005 | 0.0008 | 0.9742 | 0.9744 | 0.9714 | 0.0002 | 0.0028 | 0.9700 | 0.9701 | 0.9691 | 0.0001 | 0.0009 |
| 1 | 0.9999 | 1.0000 | 1.0000 | 0.0001 | 0.0001 | 1.0000 | 1.0000 | 1.0000 | 0.0000 | 0.0000 | 1.0000 | 1.0000 | 1.0000 | 0.0000 | 0.0000 |
| | Mean Error % | | | 0.0004 | 0.0013 | Mean Error % | | | 0.0003 | 0.0148 | Mean Error % | | | 0.0002 | 0.0052 |

Table 2. shows that the validation for standardized steady-state Concentration of oxygen $v$ along with numerical solutions for different solutions of $\gamma_{S1}$ and for values of $v_{OX} = \dfrac{-1}{2},\ \alpha = 0.01,\ \beta = 1.15 \text{ and } \gamma_{E1} = 10$

**Table 3** : **Validation of standardized steady state of Gluconic acid Concentration** $w$

| R | Gluconic acid Concentration $w$ | | | | | | | | | | | | | | |
|---|---|---|---|---|---|---|---|---|---|---|---|---|---|---|---|
| | If $\gamma_{S1}=5$ | | | | | If $\gamma_{S1}=20$ | | | | | If $\gamma_{S1}=40$ | | | | |
| | Numerical Solution Eq.29 | Eq.20 VIM [20] | Eq.29 AGM | % of error deviation of VIM | % of error deviation of AGM | Numerical Solution Eq.29 | Eq.20 VIM [20] | Eq.29 AGM | % of error deviation of VIM | % of error deviation of AGM | Numerical Solution Eq.29 | Eq.20 VIM [20] | Eq.29 AGM | % of error deviation of VIM | % of error deviation of AGM |
| 0 | 0.1820 | 0.1774 | 0.1764 | 0.0252 | 0.0307 | 0.7063 | 0.7096 | 0.7054 | 0.0046 | 0.0012 | 1.3333 | 1.4192 | 1.3992 | 0.0644 | 0.0494 |
| 0.2 | 0.1749 | 0.1705 | 0.1724 | 0.0251 | 0.0142 | 0.6790 | 0.6821 | 0.6811 | 0.0045 | 0.0030 | 1.2828 | 1.3643 | 1.3243 | 0.0635 | 0.0323 |
| 0.4 | 0.1536 | 0.1498 | 0.1400 | 0.0247 | 0.0885 | 0.5967 | 0.5993 | 0.5932 | 0.0043 | 0.0058 | 1.1305 | 1.1986 | 1.1586 | 0.0602 | 0.0248 |
| 0.6 | 0.1177 | 0.1149 | 0.1100 | 0.0237 | 0.0654 | 0.4579 | 0.4597 | 0.4532 | 0.0039 | 0.0141 | 0.8709 | 0.9195 | 0.9025 | 0.0558 | 0.0362 |
| 0.8 | 0.0668 | 0.0652 | 0.0621 | 0.0538 | 0.0683 | 0.2601 | 0.2610 | 0.2600 | 0.0034 | 0.0004 | 0.4969 | 0.5221 | 0.5021 | 0.0507 | 0.0104 |
| 1 | 0.0000 | 0.0000 | 0.0000 | 0.0000 | 0.0000 | 0.0000 | 0.0000 | 0.0000 | 0.0000 | 0.0000 | 0.0000 | 0.0000 | 0.0000 | 0.0000 | 0.0000 |
| | Mean Error % | | | 0.0254 | 0.0445 | Mean Error % | | | 0.0034 | 0.0245 | Mean Error % | | | 0.0491 | 0.1322 |

Table 3 discussed the Validation of standardized steady state concentration of gluconic acid $w$ along with numerical solutions for multiple values of $\gamma_{S1}$ and for fixed values of $\alpha=0.1,\ \beta=1\ and\ \gamma_E=5$

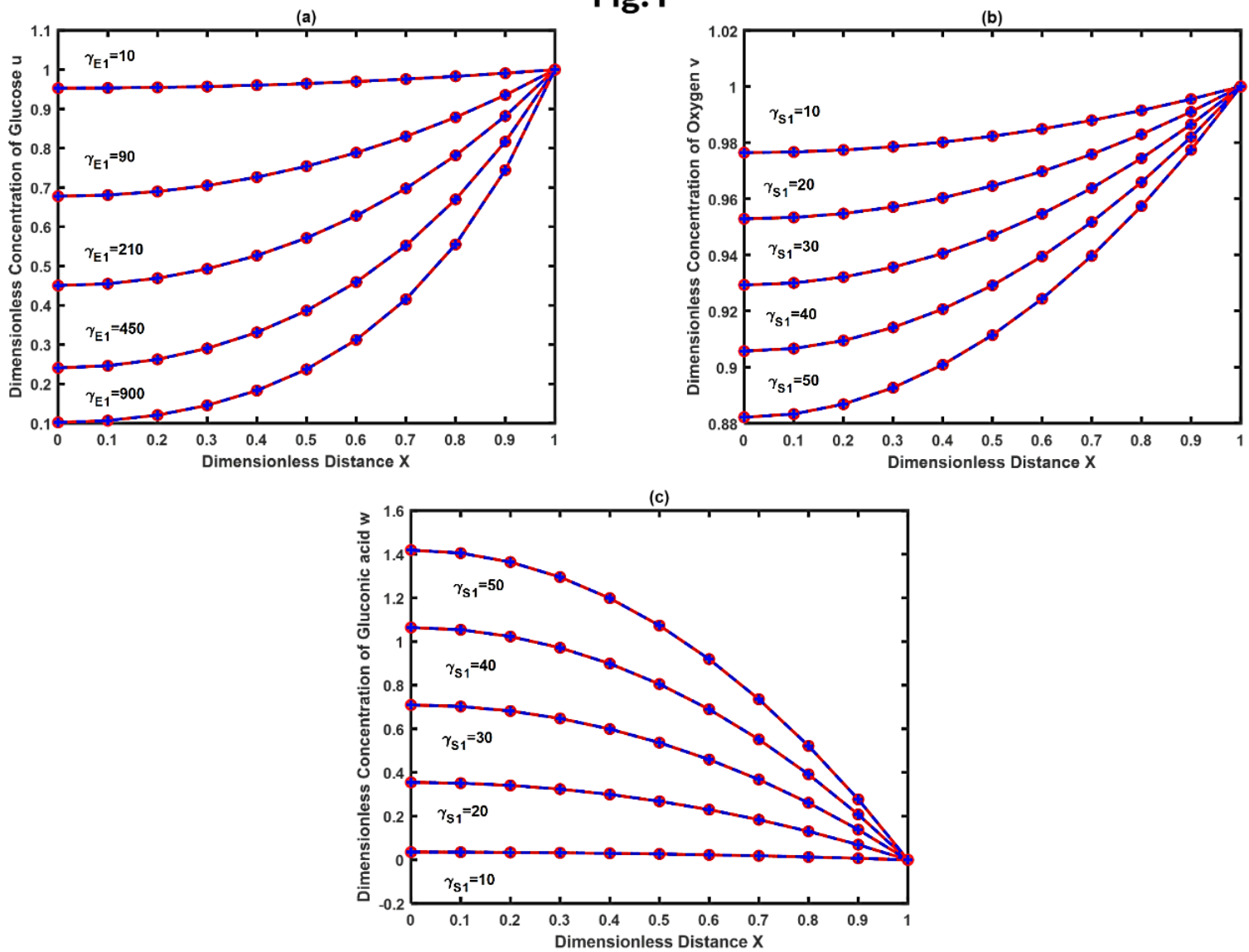

**Fig.1**

**Fig. 1** Comparison of analytic expression for the concentration of glucose , oxygen and gluconic acid and with numerical results for different parameters.

(a) $v_g = -1, \ \alpha = 0.01, \ \beta = 1.15$

(b) $v_{OX} = \dfrac{-1}{2}, \ \alpha = 0.01, \ \beta = 1.15 \ and \ \gamma_{E1} = 10$

(c) $v_a = 1, \ \alpha = 0.1, \ \beta = 1 \ and \ \gamma_{E1} = 5$ dash dotted lines , spotted lines signify the New iterative solution ,New analytical solution along with solid lines signify the numerical results.

**Fig.2**

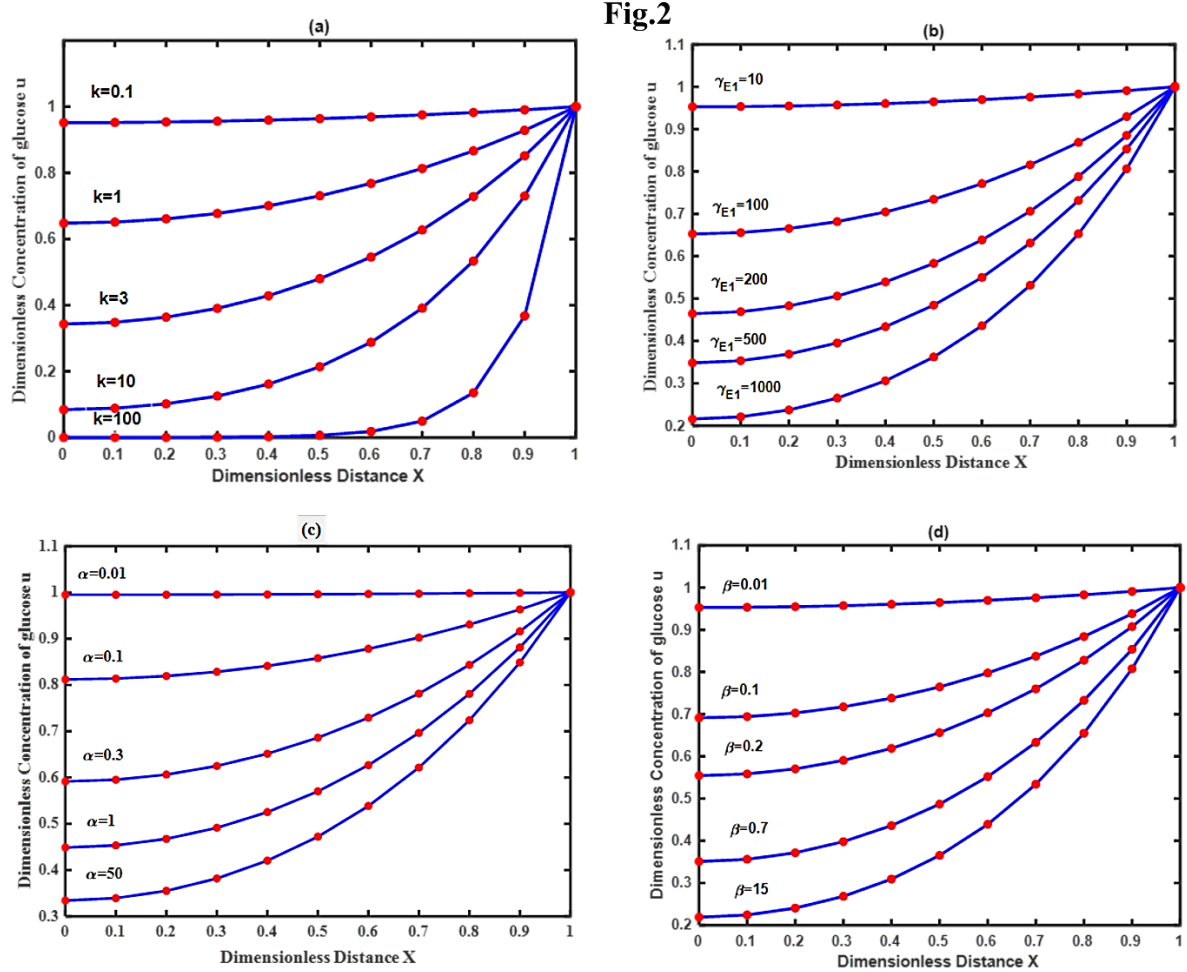

**Fig. 2** Plot for the concentration profiles for glucose *u(X)* and non-dimensional distance X calculated via Eqn. 20 for the values of:

(a) *parameter k*

$(b) \alpha = 0.01, \beta = 1.15$

$(c) \gamma_{E1} = 10, \beta = 1.15$

$(d) \gamma_{E1} = 10, \alpha = 1$

**Fig.3**

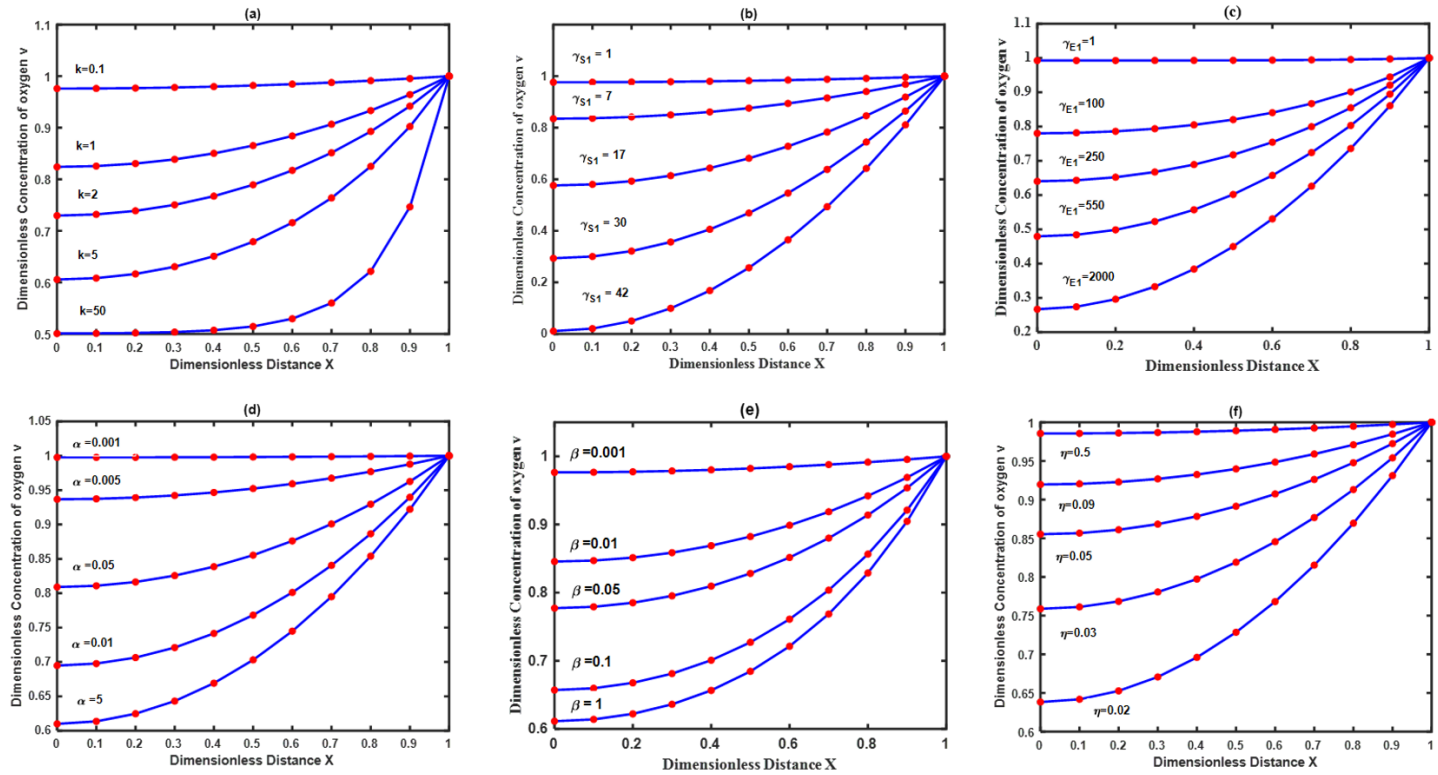

**Fig:3** Plot for the concentration profiles of oxygen $v(X)$ and non-dimensional distance X calculated via Eqn. 21 for values of:

$(a)\gamma_{S1} = 1, \gamma_{E1} = 10, \eta = 0.1$

$(b)\gamma_{E1} = 10, \eta = 0.1, \alpha = 0.01, \beta = 1$

$(c)\gamma_{S1} = 30, \eta = 0.1, \alpha = 0.01, \beta = 1$

$(d)\gamma_{S1} = 10, \gamma_{E1} = 1, \eta = 0.1, \beta = 0.01$

$(e)\gamma_{S1} = 4, \gamma_{E1} = 1, \eta = 0.1, \alpha = 0.1$

$(f)\gamma_{S1} = 3, \gamma_{E1} = 4, \alpha = 0.01, \beta = 1$

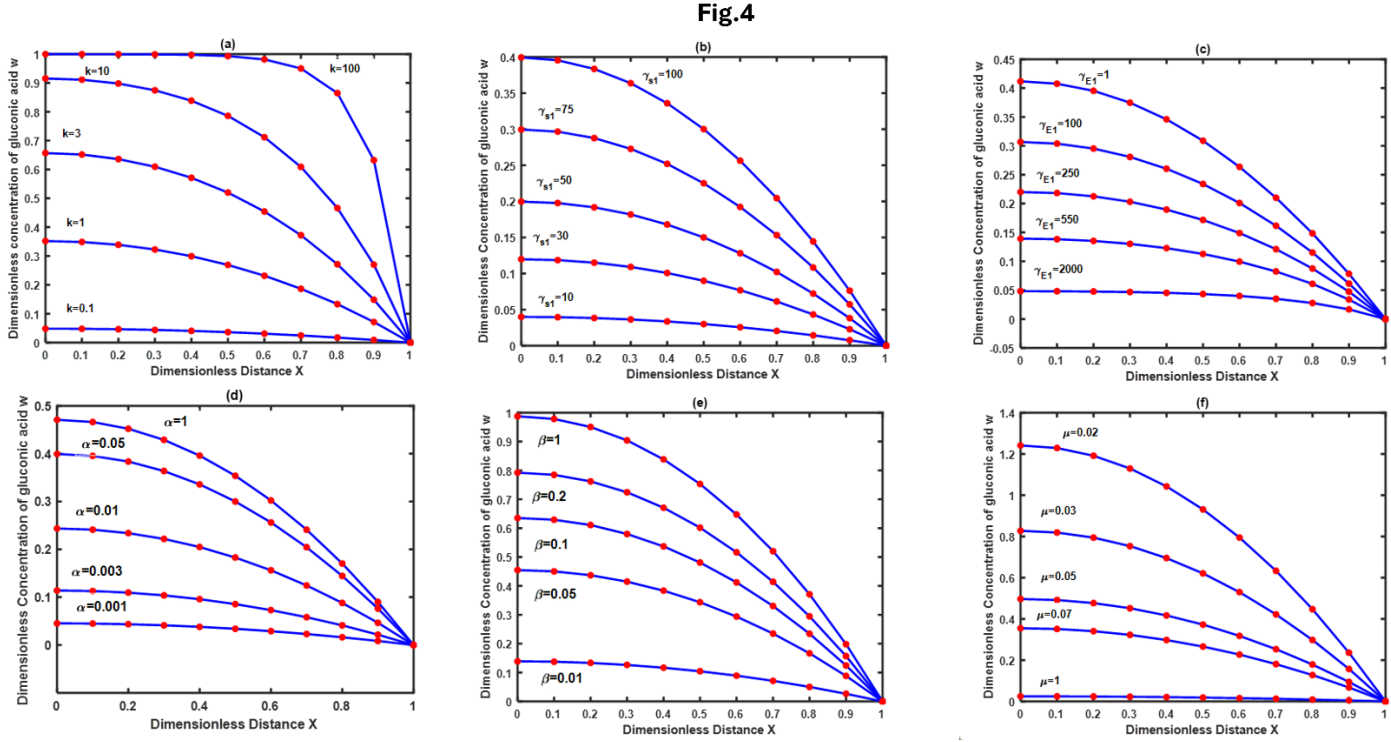

**Fig:4** Plot for the concentration profiles of gluconic acid $w(X)$ along with non-dimensional distance X calculated via Eqn. 22 for values of:

$(a)\, \gamma_{E1} = 10, \mu = 1, \gamma_S = 10$

$(b)\, \gamma_{E1} = 10, \mu = 1, \alpha = 0.05, \beta = 0.01$

$(c)\, \gamma_{S1} = 100, \mu = 1, \alpha = 0.05, \beta = 0.01$

$(d)\, \gamma_{S1} = 100, \gamma_{E1} = 10, \mu = 1, \beta = 0.01$

$(e)\, \gamma_{S1} = 32, \gamma_{E1} = 10, \mu = 1, \alpha = 0.1$

$(f)\, \gamma_{S1} = 10, \gamma_{E1} = 1, \alpha = 0.01, \beta = 0.01$

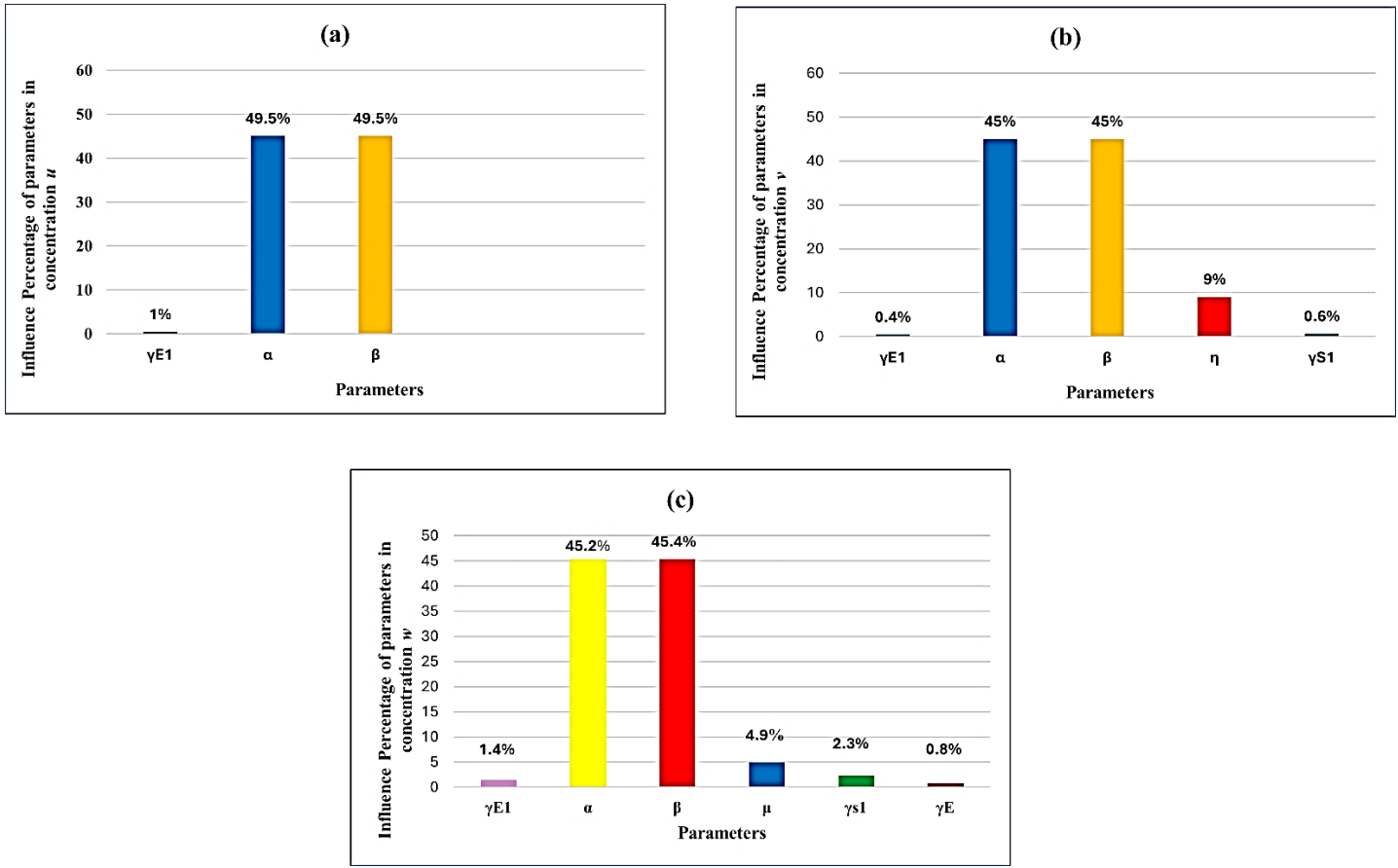

**Fig. 5** Impact  proportion of the parameters in various concentrations
*(a) u   , (b) v  ,(c) w*

# Fig.6

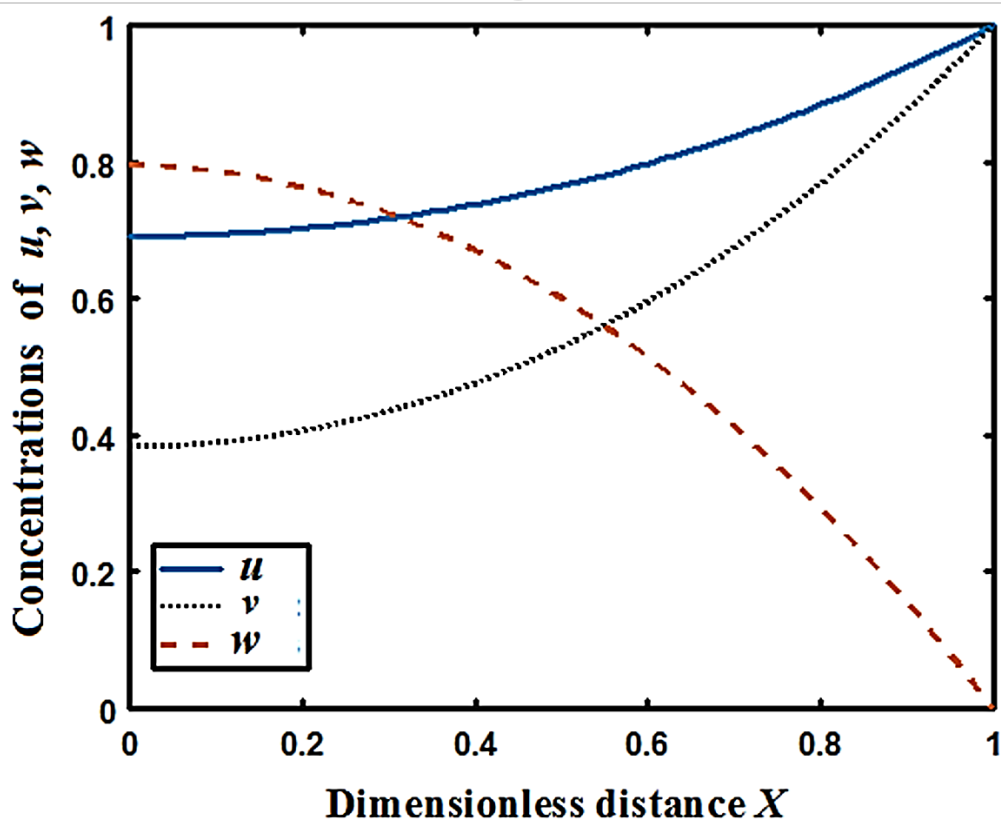

**Fig.6** Dimensionless Concentration of u, v, w versus Dimensionless Distance X and solid line, dash lines, and dotted lines represent the Analytical solution for Concentration of glucose(Eq.20), oxygen(Eq.21), and gluconic acid(Eq.22) respectively.

## Authors statements

***Availability of data and material:*** The data sets used and/or analysed during the current study are available from the corresponding author upon reasonable request.

***Funding information:*** The Deanship of Graduate Studies and Scientific Research at Jazan University in Saudi Arabia, under Project number GSSRD-24.

**Conflict of interests:** No potential conflict of interest was reported by the author(s).

**Ethics and Consent:** Not applicable.

**Consent for publication:** Informed consent was obtained from all individual participants included in the study.

**Author's contributions:** Conceptualization, K. Saranya; Methodology, M. Suguna; Software, M. Suguna; Validation, K. Saranya; Formal analysis, Salahuddin; Investigation, Salahuddin; Resources, M. Suguna; Data curation, Salahuddin; Writing—original draft preparation, M. Suguna; Writing—review and editing, Salahuddin; Visualization, Salahuddin; Supervision, K. Saranya; Project administration, Salahuddin; Funding acquisition, Salahuddin. All authors have read and agreed/approved to the published version of the manuscript.

*Acknowledgements:* *The authors gratefully acknowledge the funding of the Deanship of Graduate Studies and Scientific Research, Jazan University, Saudi Arabia, through Project Number: GSSRD-24.*

*Author's information:* *Not applicable.*